%\magnification=\magstep1
\input amstex
\documentstyle{amsppt}
%%%%%%%%%%%%%%%%%%%%%%%%%%%%%%%%
\def\pf{\hfill\hfill\qed}

\def\br{\Bbb R}

%\def\deffy{\diff_{E}(M)}
%%%%%%%%%%%%%%%%%%%%%%%%%%%%%%%%
\leftheadtext{Augustin Banyaga}
\TagsOnRight
\NoBlackBoxes
\topmatter
\title
ON THE GROUP OF STRONG SYMPLECTIC HOMEOMORHISMS
\endtitle
\vskip .5in
\author
Augustin Banyaga
\endauthor
\vskip .5in
\abstract
We generalize the "hamiltonian topology" on hamiltonian isotopies to an intrinsic "symplectic topology" on the space 
of symplectic isotopies. We use it to define the group $SSympeo(M,\omega)$ of strong symplectic homeomorphisms, 
which generalizes the group $Hameo(M,\omega)$ of hamiltonian homeomorphisms introduced by Oh and Muller. The group
$SSympeo(M,\omega)$ is arcwise connected, is contained in the identity component of $Sympeo(M,\omega)$; it
contains $Hameo(M,\omega)$ as a normal subgroup and coincides with it when $M$ is simply connected. Finally its commutator subgroup $[SSympeo(M,\omega),SSympeo(M,\omega)]$
is contained in $Hameo(M,\omega)$.
\endabstract
\endtopmatter
\document
\baselineskip 20pt

{\bf 1. Introduction}

\vskip .1in

The Eliashberg-Gromov  symplectic rigidity theorem says that the group
 $Symp(M,\omega)$ of symplectomorphisms of
a closed symplectic manifold $(M,\omega)$ is $C^0$ closed in the group $Diff^\infty(M)$ of $C^\infty$ diffeomorphisms
of $M$ ( see [8]). This means that the "symplectic" nature of a sequence of symplectomorphisms survives topological limits.
Also Lalonde-McDuff-Polterovich have shown in [9] that  for a symplectomorphism, being "hamiltonian" is topological in 
nature. These phenomenons attest that there is a {\it $C^0$  symplectic topology}  underlying the symplectic geometry
of a symplectic manifold $(M,\omega)$.
  \vskip .1in
  According to Oh-Muller ([10]), the automorphism group of the $C^0$ symplectic topology is the closure of
 the group  $Symp(M,\omega)$  in the group $Homeo(M)$ of 
 homeomorphisms of $M$ endowed with the $C^0$ topology. That group, denoted $Sympeo(M,\omega)$ has been
  called the group of symplectic
 homeomorphisms:
 $$
 Sympeo(M,\omega) = :\overline {Symp(M,\omega)}.
 $$
 \vskip .1in
The $C^0$ topology on $Homeo(M)$ coincides with the metric topology coming from the metric
$$
\overline {d}(g,h) = max ( sup_{x\in M} d_0(g(x), h(x)),  sup_{x\in M} d_0(g^{-1}(x), h^{-1}(x))
$$
where $d_0$ is a distance on $M$ induced by some riemannian metric [11].
\vskip .1in

On the space $PHomeo(M)$ of
continuous paths $\gamma : [0,1] \to Homeo(M)$, one has the distance
$$
\overline {d} (\gamma, \mu) = sup_{t\in [0,1]} \overline {d}(\gamma(t), \mu(t)).
$$

Consider the space $ PHam(M)$  of all isotopies $\Phi_H = [t\mapsto \Phi_H^t]$ 
where $\Phi_H^t$ is the family of hamiltonian diffeomorphisms obtained by integration of the family 
of vector fields $X_H$ for a smooth family $H(x,t)$ of real functions on $M$, i.e.
$$
\frac {d}{dt} \Phi_H^t(x) = X_H(\Phi_H^t(x))
$$

and $\Phi_H^0 = id$.

Recall that $X_H$ is uniquely defined by the equation
$$
i(X_H)\omega = dH
$$

where $i(.)$ is the interior product.

\vskip .1in

The set of time one maps of all  hamiltonian isotopies $\{\Phi_H^t\}$ form a group, denoted $Ham(M,\omega)$ and called
the group of hamiltonian diffeomorphisms.

\vskip .1in

{\bf Definition} {\it The hamiltonian topology} [11] on $PHam(M)$  is the metric topology
defined by the distance
$$
d_{ham} (\Phi_H,\Phi_{H'}) = ||H-H'|| + \overline {d}(\Phi_H,\Phi_{H'})
$$
where
$$ 
||H-H'||  = \int_0^1 osc (H - H') dt.
$$

and the oscillation of a function $u$ is
 $$
 osc (u) = max_{x \in M} u(x) - min_{x\in M} u(x).
 $$
 
\vskip .1in

Let $Hameo(M,\omega)$ denote the space of all homeomorphisms $h$ such that there exists a 
continuous path $\lambda \in PHomeo(M)$ such that

$ \lambda(0) = id$, $\lambda(1) = h$

 and 
there exists a Cauchy sequence (for the $d_{ham}$ distance) of hamiltonian isotopies $\Phi_{H^n}$,
which $C^0$ converges to $\lambda$ ( in the $\overline {d}$ metric).
 
The following is the first important theorem in the $C^0$ symplectic topology [11]:

\vskip .1in

{\bf Theorem (Oh-Muller)}

\vskip .1in

{\it The set $Hameo(M,\omega)$ is a topological group. It is a normal subgroup of the identity component 
 $Sympeo_0(M,\omega)$ in $Sympeo(M,\omega)$.
If $H^1(M,\br) \neq 0$, then $Hameo(M,\omega)$ is strictly contained in $Sympeo_0(M,\omega)$}.

\vskip .1in

{\bf Remark}

It is still unkown  in general if the inclusion
$$
Hameo(M,\omega) \subset Sympeo_0(M,\omega)
$$
is strict.

\vskip .1in

The group $Hameo(M,\omega)$ is the topological analogue of the group $Ham(M,\omega)$ of hamiltonian diffeomorphisms.
 \vskip .1in

The goal of this paper is to construct a subgroup of $Sympeo_0(M,\omega)$,  denoted $SSympeo(M,\omega)$
and nicknamed the group of strong symplectic homeomorphisms, containing $Hameo(M,\omega)$, that is:
$$
Hameo(M,\omega) \subset SSympeo(M,\omega) \subset Sympeo_0(M,\omega).
$$

\vskip .1in

Like $Hameo(M,\omega)$, the group $SSympeo(M,\omega)$ is defined using a blend
of the $C^0$ topology and the Hofer topology on the space $Iso(M, \omega)$ of symplectic isotopies of $(M,\omega)$. 

\vskip .1in

We believe that $SSympeo(M,\omega)$ is "more right" than the group $Sympeo(M,\omega)$ for the $C^0$ symplectic topology. In particular the flux homomorphism seems to exist on $SSympeo(M,\omega)$. This will be the object of a futur paper.

The results of this paper have been announced in [1].

\vskip .1in

The $C^0$ counter part  of the $C^\infty$ contact topology is been worked out in [5], [6].

\vskip .3in

{\bf 2. The symplectic topology} on $Iso(M,\omega)$

\vskip .1in

Let $Iso(M,\omega)$ denote the space of symplectic isotopies of a compact symplectic manifold $(M,\omega)$.
 Recall that a symplectic isotopy is a smooth map $H : M \times [0,1] \to M$ such 
that for all $t \in [0,1]$, $h_t: M \to M, ~~~~~~~~~x\mapsto H(x,t)$ is a symplectic diffeomorphism 
and $h_0 = id$.

 The "Lie algebra" of $Symp(M,\omega)$ is
 the space $symp(M, \omega)$ of symplectic vector fields, i.e the set of vector fields
 $X$ such that $i_X\omega$ is a closed form.
 
 \vskip .1in
 
  Let $\phi_t$ be a symplectic isotopy, then
$$
\dot {\phi}_t(x) = \frac{d\phi_t}{dt}(\phi^{-1}_t(x))
$$

is a smooth family of symplectic vector fields.

\vskip .1in

By the theorem of existence and uniqueness of solutions of ODE's, 
$$
 \Phi \in Iso(M, \omega) \mapsto \dot {\phi}_t
$$

is a 1-1 correspondence between $Iso(M,\omega)$ and the space $C^\infty ([0,1], symp(M,\omega))$ 
of smooth families of symplectic vector fields.
 Hence any distance on $C^\infty ([0,1], symp(M,\omega))$ gives rise to
a distance on $Iso(M,\omega)$.

\vskip .1in

{\bf An intrinsic topology on  the space of symplectic vector fields.}

\vskip .1in

 We define a norm $||.||$ on $symp(M,\omega)$ as follows: first we fix a riemannian metric $g$ (which may be 
 the one we used to define $d_0$ above, or any other riemannian metric), and a basis 
 $\Cal B$  = $\{h_1,..,h_k\}$ of harmonic 1-forms.For Hodge theory, we refer to [12].
 
  Recall that the space $harm^1(M,g)$ of harmonic 1-forms is
 a finite dimensional vector space and its dimension is the first Betti number of $M$.
  
 On  $harm^1(M,g)$, we put the following "Euclidean" norm:
  
 for $H \in harm^1(M,g)$ , $H = \sum \lambda_i h_i$, define:
  $$
  |H|_{\Cal B} = : \sum |\lambda_i|.
  $$
  \vskip .1in

This norm is equivalent to any other norm. Here we choose this one for convenience in the calculations and estimates to come later.

\vskip .1in

 Given $X \in sym(M, \omega)$, we consider the Hodge
 decomposition of $i_X\omega$ [10] : there is a unique harmonic 1-form $H_X$ and a unique function 
 $u_X$ such that
 $$
 i_X\omega  = H_X + du_X
 $$
 
 Now we define a norm $||.||$ on the the space $symp(M,\omega)$ by:
 $$
 ||X|| = |H_X|_{\Cal B}  + osc(u_X). ~~~~~~~~~~~~~~~~~~~~~~~~~~~~~\tag 1
 $$
 
 It is easy to see that this is  a norm. Let us just verify that $||X|| = 0$ implies that $X = 0$. 
 Indeed  $|H_X|_{\Cal B} = 0$ implies that
 $i_X\omega = du_X$, and $osc(u_X) = 0$ implies that $u_X$ is a constant, therefore $du_X = 0$.
 
 \vskip .1in
 {\bf Remark}
 
This norm is not invariant by $Symp(M,\omega)$.
 Hence it does not define a Finsler metric on $Symp(M,\omega)$.
 
 \vskip .1in
 
 The norm $||.||$ defined above depends of course on the riemannian metric $g$ and the basis $\Cal B$
 of harmonic 1-forms. However, we have the following:
 
 \vskip .1in
 
 {\bf Theorem 1}
 
 \vskip .1in
 
 {\it  All the norms $||.||$ defined by equation (1) using different riemannian metrics
 and different basis of harmonic 1-forms are equivalent.
 
 Hence the topology on the space $symp(M,\omega)$  of symplectic vector fields
  defined by the norm (1) is intrinsic :  it is independent of the choice 
 of the riemannian metric $g$ and of the basis $\Cal B$ of harmonic 1-forms.}

\vskip .1in

 For each symplectic isotopy $ \Phi = (\phi_t)$, consider the Hodge decomposition of $i_{(\dot {\phi}_t)}\omega$ 
$$
 i_{(\dot {\phi}_t)}\omega = \Cal H^{\Phi}_t + du^{\Phi}_t  
$$
 
where
$\Cal H^{\Phi}_t$ is a  harmonic 1-form.
\vskip .1in

We define the length $l(\Phi)$ of the isotopy $\Phi = (\phi_t)$ by:
 $$
 l(\Phi) = \int_0^1 (|\Cal H^{\Phi}_t | + osc (u^{\Phi}_t))dt = 
 \int_0^1 ||\dot {\phi}_t|| dt 
 $$

 One also writes
 $$
 \int_0^1 ||\dot {\phi}_t|| dt = 
|||\dot {\phi}_t|||.
 $$
\vskip .1in

In the expressions above, we have written $|\Cal H^{\Phi}_t |$ for  $|\Cal H^{\Phi}_t |_{\Cal B}$, where ${\Cal B}$
is a fixed basis of $harm^1(M,g)$, for a fixed riemannian metric $g$.

\vskip .1in
We define the distance  $D_0(\Phi, \Psi)$ between two symplectic isotopies $\Phi = (\phi_t)$ and $\Psi = (\psi_t)$ by:
 $$
 D_0(\Phi, \Psi) = |||\dot {\phi}_t - \dot {\psi}_t||| = :\int_0^1 (| \Cal H_t^{\Phi}
  -\Cal H_t^{\Psi}| + osc
  (u^{\Phi_t} - u^{\Psi_t})) dt.  
 $$
Denote by  $\Phi^{-1} = (\phi_t^{-1})$ and by  $\Psi^{-1} = (\psi_t^{-1})$ the inverse isotopies.

\vskip .1in

{\bf Remarks}

1. The distance $D_0(\Phi, \Psi) \neq l(\Psi^{-1}\Phi)$ unless $\Psi$ and $\Phi$ 
are hamiltonian isotopies ( see proposition 1).

2. $l(\Phi) \neq l(\Phi^{-1})$ unless $\Phi$ is hamiltonian.

\vskip .1in

In view of the remarks above, we define a more "symmetrical" distance $D$ by:
$$
 D(\Phi, \Psi) = ( D_0(\Phi, \Psi) + D_0(\Phi^{-1}, \Psi^{-1}))/2
$$
 \vskip .1in

Following [11], we define the {\it symplectic distance } on $Iso(M, \omega)$ by:
$$
 d_{symp}(\Phi, \Psi) = \overline {d}(\Phi,\Psi) + D(\Phi,\Psi).
$$

\vskip .1in

{\bf Definition}. The {\it symplectic topology} on $Iso(M, \omega)$ is the metric topology 
defined by the distance $d_{symp}$.

\vskip .1in

{\bf Theorem 1'}

\vskip .1in

{\it The symplectic topology on $Iso(M,\omega)$  natural : it is independent of all choices involved in its definition.}

\vskip .1in

We may also define another distance $D^\infty$ on $Iso(M,\omega)$ :

$$
D_0^{\infty} (\Phi, \Psi) = sup_{t \in [0,1]} (|\Cal H_t^{\Phi} -\Cal H_t^{\Psi}| + osc
  (u^{\Phi_t} - u^{\Psi_t})) 
$$
$$
D^{\infty} (\Phi, \Psi) = ((D_0^{\infty} (\Phi, \Psi) + D_0^{\infty} (\Phi^{-1}, \Psi^{-1}))/2
$$

and

$$
d^{\infty}_{symp}(\Phi,\Psi) = \overline {d} (\Phi,\Psi) + D^{\infty}(\Phi, \Psi)
$$

\vskip .1in

\proclaim
{Proposition 1}

Let $\Phi = (\phi_t), \Psi =(\psi_t)$ be two hamiltonian isotopies and $\sigma_t = (\psi_t)^{-1}\phi_t$
then
$$ 
|||\dot {\sigma}_t||| = ||| \dot {\phi}_t - \dot {\phi}_t||| = \int _0^1 osc (u^{\Phi}_t - u^{\Psi_t})dt
$$
\endproclaim

\vskip .1in

{\bf Proof}

\vskip .1in

This follows immediately from the equation
$$
\dot {\sigma}_t = ({\psi_t}^{-1})_* (\dot {\phi}_t - \dot {\phi}_t),
$$
which is a consequence of proposition 4.
\pf

\proclaim
{Corollary}

The distance $d_{sym}$ reduces to the hamiltonian distance $d_{ham}$
when $\Phi$ and $\Psi$ are hamiltonian isotopies.

\endproclaim

 The {\it symplectic topology} reduces to the "hamiltonian topology" of [11]
on paths in $Ham(M,\omega)$.

\vskip .1in

{\bf A Hofer-like metric on $Symp(M,\omega)$}

\vskip .1in

For any $\phi \in Symp(M,\omega)$, define:
$$
e_0(\phi) = inf  (l(\Phi))
$$

where the infimum is taken over all symplectic isotopies $\Phi$ from  $\phi$ to the identity.
 The following result was proved in [2].

\vskip .1in

\proclaim
{Theorem}

 The map $ e : Symp(M,\omega) \to \br \cup \{\infty\}:$
 $$
  e(\phi) =: (e_0(\phi) + e_0(\phi^{-1}))/2
 $$
   
{\it is a metric on the identity component $Symp(M,\omega)_0$ in the group $Symp(M,\omega)$,
 i.e. it satisfies (i) $e(\phi) \geq 0$ and $e(\phi)= 0$ iff $\phi$ is the identity.

(ii) $e(\phi) = e((\phi)^{-1})$

(iii) $e\phi .\psi) \leq (e\phi) + e(\psi)$.

The restriction to $Ham(M,\omega)$
is bounded from above by the Hofer norm}.

\endproclaim

\vskip .1in

Recall that the Hofer norm [8]  of a hamiltonian diffeomorphism $\phi$ is
$$
||\phi||_H = inf (l(\Phi_H))
$$
where the infimum is taken over all hamiltonian isotopies $\Phi_H$ from $\phi$ to the identity.

\vskip .1in

The Hofer-like metric above depends on the choice of a riemannian metric $g$ and a basis $\Cal B$ of harmonic 1-forms. Hence it is not "natural".
However, by theorem 1, all the metrics constructed that way are equivalent; so they define a
 natural topology on $Symp(M,\omega)$.

\vskip .3in

3. {\bf Strong symplectic homeomorphisms}

\vskip .1in

{\bf Definition} : {\it A homeomorphism $h$ is said to be a strong symplectic homeomorphism if
there exists a continuous path
 $\lambda : [0,1] \to Homeo(M)$ such that $\lambda(0) = id ; 
 \lambda(1) = h$
and a
sequence $\Phi^n =(\phi^n_t)$ of symplectic isotopies, which converges  
 to $\lambda$ in the
 $C^0$ topology ( induced by the norm $\overline {d}$) and such that $\Phi^n$ 
  is Cauchy for the metric $d_{symp}$.}
   
\vskip .1in

We will denote by $SSympeo(M,\omega)$ the set of all strong symplectic homeomorphisms. This set is well defined
independently of any riemannian metric or any basis of harmonic 1-forms.

\vskip .1in

Clearly, if $M$ is simply connected, the set $SSympeo(M,\omega)$ coincides with the group $Hameo(M,\omega)$.

\vskip .1in

We denote by $SSympeo(M, \omega)^\infty$ the set defined like in $SSympeo(M,\omega)$ but replacing
the norm $d_{symp}$ by the norm $d_{symp}^{\infty}$.

\vskip .1in

Let $\Cal {P} Homeo(M)$ be the set of continuous paths $\gamma :[0,1] \to Homeo(M)$ such that $\gamma(0) = id$, and let $ \Cal {P}^\infty(Harm^1(M)$  be the space of smooth paths of harmonic 1-forms.

We have the following maps:

\vskip.1in

$A_1: Iso(M,\omega) \to \Cal {P} Homeo(M),  \Phi \mapsto \Phi(t)$

\vskip.1in
$A_2 : Iso(M, \omega) \to \Cal {P}^\infty(Harm^1(M), \Phi \mapsto \Cal H^{\Phi}_t$

\vskip .1in

$A_3 :  Iso(M,\omega) \to C^{\infty}(M \times [0,1],\br),\Phi \mapsto u^{\Phi}$

\vskip .1in

Let $\Cal Q$ be the image of the mapping $A = A_1 \times A_2 \times  A_3$ and $\overline {\Cal Q}$ the closure of $\Cal Q$ inside $\Cal {I}(M,\omega) = : \Cal {P} Homeo(M) \times \Cal {P}^\infty(Harm^1(M) \times C^{\infty}(M \times [0,1], \br)$, with the symplectic topology, which is the $C^0$ topology on the first factor and the metric topology from $D$ on the second and third factor.

\vskip .1in

 Then $SSympeo(M,\omega)$ is just the image of the evaluation map of the path at t= 1 of the image of the projection of $\Cal Q$ on the first factor. This defines a surjective map:
$$
a : \Cal Q \to SSympeo(M,\omega)
$$
The symplectic topology on $SSympeo(M,\omega)$ is the quotient topology induced by $a$.
\vskip .1in

Our main result is the following
 
\vskip .1in

\proclaim
{Theorem 2}

Let $(M,\omega)$ be a closed symplectic manifold. Then $SSympeo(M,\omega)$  
is an arcwise connected  topological group, containing $Hameo(M,\omega)$ as 
a normal subgroup, and contained in the identity 
component $Sympeo_0(M,\omega)$ of $Sympeo(M,\omega)$.

If $M$ is simply connected, $SSympeo(M,\omega) = Hameo(M,\omega)$. Finally, the commutator subgroup 
$[SSympeo(M,\omega), SSympeo(M,\omega)]$ of $SSympeo(M,\omega)$ is contained in $Hameo(M,\omega)$.

\endproclaim

\vskip .1in
 
{\bf Conjectures}

\vskip .1in

1. Let $(M,\omega)$ be a closed symplectic manifold, then

 $[SSympeo(M,\omega), SSympeo(M,\omega)] = Hameo(M,\omega)$.

\vskip .1in
2. The inclusion $SSympeo(M,\omega) \subset Sympeo_0(M,\omega)$ is strict.

\vskip .1in

3. The results in theorem 2 hold for $SSympeo(M,\omega)^\infty$.

\vskip .1in

Conjecture 3 is supported by a result of Muller asserting that $Hameo(M,\omega)$ coincides with 
$Hameo(M,\omega)^{\infty}$ which is defined by replacing the $L^{(1, \infty)}$ Hofer norm by the $L^{\infty}$ norm [8].

\vskip .1in

{\bf Measure preserving homeomorphisms}

\vskip .1in

On a symplectic $2n$ dimensional manifold $(M,\omega)$, we consider the measure 
$\mu_{\omega}$ defined by the Liouville volume $\omega^n$. Let $Homeo_0^{\mu_{\omega}}(M)$ be the identity component
in the group of homeomorphisms preserving $\mu_{\omega}$. We have:
$$
 Sympeo_0(M,\omega) \subset Homeo_0^{\mu_{\omega}}(M).
$$
Oh and Muller [11] have observed that $Hameo(M,\omega)$ is a sub-group of the kernel
 of Fathi's mass-flow homomorphism
 [7]. This is a homomorphism $\theta : Homeo_0^{\mu_{\omega}}(M) \to H_1(M,\br)/\Gamma$, 
 where $\Gamma$ is some 
 sub-group of $H_1(M,\br)$.
 Fathi proved that if the dimension of $M$ is bigger than 2, then $Ker \theta$ is a simple group.
 This leaves open the following question [11]:
 
 {\it Is $Homeo_0^{\mu_{\omega}}(S^2) = Sympoe_0(S^2,\omega)$ a simple group?}
 
 \vskip .1in
 
 But $Sympoe_0(S^2,\omega)$ contains $Hameo(S^2, \omega)$ as a normal subgroup.
 The question is to decide if the inclusion 
 $$
  Hameo(S^2, \omega) \subset Sympoe_0(S^2,\omega)
 $$ 
 is strict.  Since $SSympeo(S^2,\omega) = Hameo(S^2, \omega)$, our conjecture 2 implies that
 $Homeo_0^{\mu_{\omega}}(S^2) = Sympoe_0(S^2,\omega)$ is not a simple group, a conjecture of [9].

\vskip .1in

{\bf Questions}

\vskip .1in

1.  Is $SSympeo(M,\omega)$ a normal subgroup of $Sympeo_0(M,\omega)$?

2. Is $[Sympeo_0(M,\omega),Sympeo_0(M,\omega)]$ contained in $Hameo(M,\omega)$?

\vskip .3in

{\bf 4. Proofs of the results}

\vskip .1in

{\bf 4.1. Proof of theorem 1}

\vskip .2in

If $\Cal B$ and $\Cal B'$ are two basis of $harm^1(M,g)$, then elementary linear algebra shows that $|.|_{\Cal B}$ and
$|.|_{\Cal B'}$ are equivalent. This implies that  
 the corresponding norms on $symp(M,\omega)$ are also equivalent.

\vskip .1in

Let us now start our construction with a riemannian metric $g$ and a basis $\Cal B = (h_1,..h_k)$ of $harm^1(M,g)$.
We saw that for any $X \in symp(M,\omega)$, 
$$
i_X\omega = H_X + du_X
$$
and we wrote $H_X = \sum \lambda_i h_i$.

\vskip .1in

Let $g'$ be another riemannian metric. The $g'$-Hodge decomposition of $i_X\omega$ is:
$$
i_X\omega = H'_X + du'_X
$$
where $H'_X$ is $g'$-harmonic. 

\vskip .1in

Consider the $g'$- 
Hodge decompositions of the members $h_i$ of the basis $\Cal B$ i.e.
$$
h_i = h'_i + dv_i
$$
where $h'_i$ is $g'$ harmonic.

\vskip .1in

$\Cal B' = (h'_1,..h'_k)$ is a basis of $harm^1(M,g')$. Indeed if $\sum r_ih'_i = 0$, then
 $\sum r_ih_i = d(\sum r_iv_i)$. Hence  $\sum r_ih_i $ is identically zero because it is an exact harmonic form.
 Therefore all $r_i$ are zero since $\{h_i\}$ form a basis.
 
\vskip .1in
The 1-form
$$
H''_X = : \sum \lambda_i h'_i
$$
is a $g'$- harmonic form representing the cohomology class of $i_X\omega$. By uniqueness, $H'_X = H''_X$.
\vskip .1in

Hence 
$$
|H'_X|_{\Cal B'} = \sum |\lambda_i| = |H_X|_{\Cal B}
$$

\vskip .1in

Furthermore $H'_X = \sum \lambda_i (h_i - dv_i)  = H_X + dv$ where $v = - \sum \lambda_iv_i$.
Hence
$$
i_X\omega = H'_X + du'_X = H_X + d(v + u'_X)
$$
By uniqueness in the $g$-Hodge decomposition of $i_X\omega$,
$$
u_X =v + u'_X.
$$
Denote by $||X||_{g'}$, resp. $||X||_{g}$, the norm of $X$ using the riemannian metric $g'$ and 
the basis $\Cal B'$, resp. using the riemannian metric $g$ and the basis $\Cal B$.
Then:
$$
||X||_{g'} = |H'_X|_{\Cal B'} + osc (u'_X) = |H'_X|_{\Cal B'} + osc (u_X - v) 
$$
$$
 \leq |H'_X|_{\Cal B'} + osc (u_X) + osc(-v)
$$
$$
 = |H_X|_{\Cal B} +  osc (u_X)+ osc(v) = ||X||_{g} + osc (v).
$$

Similarly,
$$
||X||_{g} = |H_X|_{\Cal B} + osc (u_X) = |H_X|_{\Cal B} + osc (v + u'_X)
$$
$$
 \leq (|H_X|_{\Cal B} + osc ( u'_X) )+ osc (v) = ||X||_{g'} + osc (v).
$$

Setting $a =||X||_{g},  b =||X||_{g'},  c =osc (v)$, we proved $a\leq b+c$ and $b \leq a+c$. Substracting 
these inequalities, we get $a-b \leq b-a$ and $b-a \leq a- b$. This gives $a \leq b$ and $b \leq a$, i.e $a = b$.
 
\vskip .1in

We proved that given the couple $(g,\Cal B)$ of a riemannian metric $g$ and a basis of $g$-harmonic 1-forms ,
and any other riemannian metric $g'$, there is a basis $\Cal B'$ of $g'$-harmonic 1-forms so that
 $||X||_{g} =|X||_{g'}$, hence the norm $||.||$ is independent of the riemanian metric up to the equivalence relation
 due to change of basis. In conclusion, all the norms on $symp(M,\omega)$  given by formula (1) are equivalent.\pf
 
\vskip .1in

For the purpose of the proof of the main theorem, we fix a riemannian metric $g$ and a basis $\Cal B = (h_1,..,h_k)$
of $harm^1(M,g)$. The norm of a harmonic 1-form $H$ will be simply denoted $|H|$ and the norm of a symplectic vector field $X$
will be simply denoted $||X||$.

\vskip .3in

{\bf 4.2. Proof of theorem 2}

\vskip .1in

Let $h_i \in SSympeo(M,\omega)$ $i =1,2$ and let  $\lambda_i$  be continuous paths in $Homeo(M)$
 with $\lambda_i(0) = id$, $ \lambda_i (1) = h_i$ and let $\Phi^n_i$ be $d_{symp}$ - Cauchy sequences of symplectic isotopies 
 converging $C^0$ to $\lambda_i$.
 
Then $\Phi^n_1. (\Phi^n_2)^{-1}$ converges $C^0$ to the path $\lambda_1(t)(\lambda_2(t))^{-1}$.
Here $\Phi^n_1. (\Phi^n_2)^{-1}(t) =\phi_1^n(t) .(\phi_2^n(t))^{-1}$.

\vskip .1in

By definition of the distance $d_{symp}$, $\Phi^n$ is a $d_{symp}$ - Cauchy sequence if and only if 
 both $\Phi^n$ and $(\Phi^n)^{-1}$ are $D_0$ - Cauchy and $\overline {d}$- Cauchy sequences.

\proclaim
{Main lemma}

 If  $\Phi^n = (\phi^n_t)$ and  $\Psi^n_t = (\psi^n_t)$ are  $d_{symp}$ - Cauchy sequences in $Iso(M)$, 
 so is $\rho^n_t = \phi^n_t\psi^n_t$.
\endproclaim

\vskip .1in

It will be enough to prove that   $\rho^n_t$ is a  $D_0$ - Cauchy sequence. Indeed since $(\Phi^n)^{-1}$
 and $(\Psi^n)^{-1}$ are $D_0$ - Cauchy by assumption, the main lemma applied to their product implies that
  their product  is also $D_0$ Cauchy.
Hence $(\Psi^n)^{-1}(\Phi^n)^{-1} = (\Phi^n\Psi^n)^{-1} = (\rho^n_t)^{-1}$ is a $D_0$ - Cauchy sequence.
 This will conclude the proof
that $SSympeo(M,\omega)$ is a group.

\vskip .1in

We will use the following estimate:

\vskip .1in

 \proclaim
 {Proposition 2}
 There exists a constant $E$ 
 such that for any  $ X\in symp(M,\omega)$, and $H \in harm^1(M,g)$
 $$
 |H(X)| =: sup _{x\in M}  |H(x)(X(x))| \leq E||X||. |H|
 $$
 \endproclaim
 \vskip .1in
 
 \demo
 {Proof}
 Let $(h_1,..,h_r)$ be the chosen basis for harmonic 1-forms and let  $E = max_i E_i$
 and $E_i =  sup_{V} (sup_{x \in M} |h_i(x)(V(x))|$ where $V$ runs over all symplectic
 vector fields $V$ such that $||V|| = 1$.
 \vskip .1in
 Without loss of generality, we may suppose $X \neq 0$ and set $V = X/||X||$. Let $H = \sum \lambda_i h_i$. Then
$H(X) = ||X|| \sum \lambda_i h_i(V)$.  Hence 
$$
|H(X)| \leq ||X|| \sum |\lambda_i| sup_x (|h_i(x)(V)(x)|)
\leq ||X|| \sum |\lambda_i|E = E||X||.|H|.
$$

\enddemo

\pf

\vskip .1in

We will also need the following standard facts: 

\proclaim
{Proposition 3}

Let $\phi$ be a diffeomorphism, $X$  a vector field and $\theta$ a differential form on a smooth manifold $M$,
Then
$$
(\phi^{-1})^*[ i_X \phi^*\theta] = i_{\phi_*X} \theta
$$
\endproclaim

\vskip .1in 
\proclaim
{Proposition 4}

If $\phi_t, \psi_t$ are any isotopies, and
 if we denote by $\rho_t = \phi_t\psi_t$, and by $\underline {\phi}_t = (\phi)^{-1}_t$  then
$$
\dot {\rho}_t = \dot {\phi}_t + (\phi_t)_* \dot {\psi}_t
$$
and
$$
\dot {\underline {\phi}}_t = -((\phi)^{-1}_t)_* (\dot {\phi}_t)
$$
\endproclaim

\vskip .1in

\proclaim
{Proposition 5}

Let $\theta_t$ be a smooth family of closed 1-forms and $\phi_t$ an isotopy,
then
$$
\phi_t^*\theta_t - \theta_t = dv_t
$$
where
$$
v_t = \int_0^t (\theta_t(\dot {\phi}_s) \circ \phi_s)ds
$$
\endproclaim
\vskip .3in

{\bf Proof of the main lemma}

\vskip .2in

If $\phi_t, \psi_t$ are symplectic isotopies, and if $\rho_t = \phi_t \psi_t$, propositions 3, 4 and 5 give: 
 
$$
i(\dot {\rho}_t)\omega = \Cal H^{\Phi}_t +\Cal H^{\Psi}_t + d K(\Phi, \Psi)
$$

where $K = K(\Phi,\Psi) = u^{\Phi}_t + (u^{\Psi}_t) \circ (\phi_t)^{-1} + v_t(\Phi, \Psi)$, and

$$
v_t(\Phi, \Psi) = \int_0^t (\Cal H^{\Psi}_t (\dot {\underline{\phi}}_s) \circ  \phi_s^{-1}) ds.
$$

Let now $\phi_t^n, \psi_t^n$ be Cauchy sequences of symplectic isotopies, 
and consider the sequence $\rho^n_t =\phi_t^n \psi_t^n$.

We have:
$$
||| \dot {\rho}^n_t -\dot {\rho}^m_t||| =  \int_0^1 |\Cal H^{{\Phi}^n}_t 
- \Cal H^{{\Phi}^m}_t + \Cal H^{{\Psi}^n}_t - \Cal H^{{\Psi}^m}_t | + osc(K(\Phi^n,\Psi^n) -K(\Phi^m,\Psi^m) )dt
$$

$$
\leq \int_0^1 |\Cal H^{{\Phi}^n}_t - \Cal H^{{\Phi}^m}_t)|dt +\int_0^1 
|\Cal H^{{\Psi}^n}_t - \Cal H^{{\Psi}^m}_t)|dt
$$

$$
 + \int_0^1 osc(u^{\Phi^n}_t -u^{\Phi^m}_t)dt +\int_0^1 
osc(u^{\Psi^n}_t) \circ (\phi^n_t)^{-1} -u^{\Psi^m}_t \circ (\phi^m_t)^{-1}) dt
$$

$$
 + \int_0^1 osc( v_t(\Phi^n, \Psi^n) - v_t(\Phi^m, \Psi^m)dt
$$

$$
= |||\dot {{\phi}^n}_t -\dot {\phi^m}_t|||+\int_0^1 |\Cal H^{\Psi^n}_t -
 \Cal H^{\Psi^m}_t)|dt + A + B
$$
where 
$$
A = \int_0^1 osc(u^{\Psi^n}_t) \circ (\phi^n_t)^{-1} -u^{\Psi^m}_t 
\circ (\phi^m_t)^{-1}) dt
$$ 
and
$$
B =\int_0^1 osc(v_t(\Phi^n, \Psi^n)  - v_t(\Phi^m, \Psi^m)dt
$$ 
\vskip .1in
We have:
$$
A \leq \int_0^1 osc(u^{\Psi^n}_t) \circ (\phi^n_t)^{-1} -u^{\Psi^m}_t \circ (\phi^n_t)^{-1}) dt +
\int_0^1 osc(u^{\Psi^m}_t) \circ (\phi^n_t)^{-1} - (u^{\Psi^m}_t) \circ (\phi^m_t)^{-1}) dt
$$
$$
= \int_0^1 osc(u^{\Psi^n}_t -u^{\Psi^m}_t) dt + C
$$
where
$$
C = \int_0^1 osc(u^{\Psi^m}_t \circ (\phi^n_t)^{-1} -u^{\Psi^m}_t
 \circ (\phi^m_t)^{-1}) dt.
$$

\vskip .1in

 Hence
 $$
||| \dot {\rho}^n_t -\dot {\rho}^m_t||| \leq |||\dot {\phi}^n_t -\dot {\phi}^m_t|||
$$

$$
 +\int_0^1 |\Cal H^{\Psi^n}_t - \Cal H^{\Psi^m}_t)|dt + \int_0^t osc(u^{\Psi^n}_t -u^{\Psi^m}_t) dt  + B + C
$$

$$
=|||\dot {\phi}^n_t -\dot {\phi}^m_t||| +|||\dot {\psi}^n_t -\dot {\psi}^m_t|||+ B + C
$$

 \vskip .1in
 
We now show that $C \to 0$ when $m,n \to \infty$.

\vskip .1in

{\bf Sub-lemma 1 (reparametrization lemma [11])}

\vskip .1in

{\it $\forall \epsilon \ge 0, \exists m_0$ such that
$$
 C =\int_0^1 osc(u^{\Psi^m}_t \circ (\phi^n_t)^{-1} -u^{\Psi^m}_t
 \circ (\phi^m_t)^{-1}) dt =: ||u^{\Psi^m}_t \circ (\phi^n_t)^{-1} - u^{\Psi^m}_t\circ (\phi^m_t)^{-1})||
 \leq \epsilon
 $$ 
 
if $m \geq m_0$ and $n$ large enough}
  
\vskip .1in

{\bf Remark}

\vskip .1in

This is the "reparametrization lemma" of Oh-Muller [11] (lemma 3.21. (2)). For the convenience of the reader
and further references, we include their proof.

{\bf Proof}

\vskip .1in

For short, we write $u_m$ for $u^{\Psi^m}_t$ and $\mu^n_t$ for $(\phi^n_t)^{-1}$.

First, there exists $m_0$ large such that $||u_m - u_{m_0}|| \leq \epsilon/3$  for $m\geq m_0$, 
since  $(u_m)$ is a Cauchy  sequence for the distance  $d(u_n, u_m)  = \int_0^1  osc (u_n - u_m)dt$.

Therefore 
$$
||u_m\circ \mu^n_t - u_m\circ \mu^m_t))|| \leq ||u_m\circ \mu^n_t
  - u_{m_0}\circ \mu^n_t))|| + ||u_{m_0}\circ \mu^n_t
  - u_{m_0}\circ \mu^m_t))|| +||u_{m_0}\circ \mu^m_t - u_m\circ \mu^m_t))|| 
$$
$$
= ||u_m - u_{m_0} || + ||u_{m_0}\circ \mu^n_t - u_{m_0}\circ \mu^m_t))|| + ||u_{m_0} - u_m||
$$
$$
\leq (2/3)\epsilon + ||u_{m_0}\circ \mu^n_t - u_{m_0}\circ \mu^m_t))|.
$$
 
By uniform continuity of $u_{m_0}$, there exists a  positive $\delta$ such that if
 $\overline {d}(\mu^m_t,\mu^n_t) \leq \delta$, then
  max $osc ((u_{m_0}\circ \mu^n_t - u_{m_0}\circ \mu^m_t)) \leq \epsilon/6$.
  Hence $||u_{m_0}\circ \mu^n_t - u_{m_0}\circ \mu^m_t))|| \leq \epsilon/3$ for $n,m$ large.
  Recall that $\mu^n_t$ is a $\overline {d}$- Cauchy sequence. \pf

\vskip .1in

To show that  $\dot {\rho}^n_t$ is a Cauchy sequenence, the only thing which is left
is to show that $B \to 0$ when $n,m\to \infty$.

\vskip .1in

Let us denote $v_t(\Phi^n, \Psi^n)$ by $v_t^n$ , $\Cal H^{\Psi^n}_t$  by $\Cal H_n^t$  or $\Cal H_n$ 
and $(\phi^n_t)^{-1}$ by $\mu^n_t$.

\vskip .1in

For a function on $M$, we consider the norm
$$
|f| = sup_{x \in M} |f(x)|
$$

\vskip .1in

We have:
 $$
 |v^n_t - v^m_t| =| \int_0^t (\Cal H_n(\dot {\mu}_s^n)\circ \mu_s^n -
 \Cal H_m(\dot {\mu}_s^m)\circ \phi_s^m)ds|
 $$
 
 $$
\leq  \int_0^1 |((\Cal H_n - \Cal H_m)(\dot {\mu}_s^n))\circ \mu^n_s |ds
$$

$$
+\int_0^1 |\Cal H_m ( \dot {\mu}_s^n - \dot {\mu}_s^m))\circ \mu_s^m|ds
$$

$$
 + \int_0^1 |\Cal H_m( \dot {\mu}_s^n) \circ \mu_s^n - \Cal H_m( \dot {\mu}_s^n )\circ \mu_s^m|ds
$$

The last integral can be estimated as follows:

$$
\int_0^1 |\Cal H_m( \dot {\mu}_s^n) \circ \mu_s^n - \Cal H_m( \dot {\mu}_s^n )\circ \mu_s^m|ds
$$

$$
\leq \int_0^1 |\Cal H_m( \dot {\mu}_s^n) \circ \mu_s^n - \Cal H_m( \dot {\mu}_s^{n_0} )\circ \mu_s^n|ds ~~~~~~~~~~~\tag 1
$$

$$
+ \int_0^1 |\Cal H_m( \dot {\mu}_s^{n_0}) \circ \mu_s^n - \Cal H_m( \dot {\mu}_s^{n_0} )\circ \mu_s^m|ds ~~~~~~~\tag 2
$$

$$
+ \int_0^1 |\Cal H_m( \dot {\mu}_s^{n_0}) \circ \mu_s^m- \Cal H_m( \dot {\mu}_s^n )\circ \mu_s^m|ds ~~~~~~~~~~\tag 3
$$

for some integer $n_0$.

\vskip .1in

Proposition 2 gives  $E|\Cal H_m| D_0((\Phi^n)^{-1}, (\Phi^{n_0})^{-1}) \leq
 2E|\Cal H_m| D((\Phi^n), (\Phi^{n_0})^{-1})$ as an upper bound for (1) and (3).

\vskip .1in

It also gives the following estimates:
 
\vskip .1in
$$
\int_0^1 |((\Cal H_n - \Cal H_m)(\dot {\mu}_s^n))\circ \mu^n_s | ds
\leq E |\Cal H_n - \Cal H_m|\int_0^1 || \dot {\mu}_s^n)|| ds
$$

$$ 
= E.|\Cal H_n - \Cal H_m|.l((\Phi^n)^{-1})
$$

and
 
$$
\int_0^1 |(\Cal H_m( \dot {\mu}_s^n - \dot {\mu}_s^m))\circ \mu_s^m|ds 
\leq E.|\Cal H_m| \int_0^1 ||( \dot {\mu}_s^n - \dot {\mu}_s^m)|| ds
$$

$$
= E|\Cal H_m|D_0((\Phi^n)^{-1} ,(\Phi^m))^{-1}) \leq 2E|\Cal H_m|D(\Phi^n ,\Phi^m).
$$

Therefore, we get the following estimate:
$$
|v^n_t - v^m_t| \leq E.|\Cal H_n - \Cal H_m| l(\Phi^n)^{-1}) +
 E|\Cal H_m|2(D(\Phi^n ,\Phi^m) + 4 D(\Phi^n, \Phi^{n_0})) +  G
$$

where 
$$
 G = \int_0^1 |\Cal H_m( \dot {\mu}_s^{n_0}) \circ \mu_s^n - \Cal H_m( \dot {\mu}_s^{n_0} )\circ \mu_s^m|ds
$$

\vskip .1in

Since $osc(v^n_t - v^m_t) \leq 2|v^n_t - v^m_t|$, we see that

$$
\int_0^1 osc(v^n_t - v^m_t)dt \leq 2(l(\Phi^n)^{-1}) \int_0^1 |\Cal H^t_n - \Cal H^t_m|dt 
$$

$$
+ E2(D(\Phi^m,\Phi^n) + 4 D(\Phi^n,\Phi^{n_0}) \int_0^1 |\Cal H^t_m| dt) + \int_0^1 G dt
$$

\vskip .1in
We need the following facts:

\vskip .1in

{\bf Sub-lemma 2 (Reparametrization lemma)}

\vskip .1in

{\it $\forall \epsilon \ge 0, \exists n_0$ such that
$$
 L = \int_0^1 G dt =\int_0^1 (\int_0^1 |\Cal H_m( \dot {\mu}_s^{n_0}) \circ \mu_s^n - \Cal H_m( \dot {\mu}_s^{n_0} )
 \circ \mu_s^m|ds) dt \leq \epsilon
$$
 
for  $n \geq n_0$ and $m$ sufficiently large.}

\vskip .1in

{\bf Proposition 6}

\vskip .1in

{$l((\Phi^n))^{-1}$  and  $\int_0^1 |\Cal H^t_m| dt$  are bounded for every $n,m$.}

\vskip .1in

We finish first the estimate for $\int_0^1 osc(v^n_t - v^m_t)dt$ using  sub-lemma 2 and proposition 6.

\vskip .1in

Putting together all the information we gathered, we see that:

$$
\int_0^1 osc(v^n_t - v^m_t)dt \leq 2(l(\Phi^n)^{-1}) \int_0^1 |\Cal H^t_n - \Cal H^t_m|dt 
$$

$$
+ E(2D(\Phi^m,\Phi^n) + 4 D(\Phi^n,\Phi^{n_0})( \int_0^1 |\Cal H^t_m| dt) + L
$$

$$
\leq 2 l((\Phi^n)^{-1}) D(\Phi^n, \Phi^m)+ E(2D(\Phi^m,\Phi^n) +4D(\Phi^n,\Phi^{n_0}) \int_0^1 |\Cal H^t_m| dt + L
$$

Therefore:
$$
\int_0^1 osc(v^n_t - v^m_t)dt \to 0
$$ 

when $n, m \to \infty$, and $n_0$  is chosen sufficiently large. This finishes the proof of the main lemma.\pf

\vskip .1in

{\bf Proof of proposition 6}

\vskip .1in

This follows from the estmates:
 $$
 l((\Phi^n)^{-1}) \leq D((\Phi^n)^{-1}, \Phi^{n_0}) + l(\Phi^{n_0})
 $$
 and
 
 $$
 \int_0^1 |\Cal H^t_m| dt \leq \int_0^1 |\Cal H^t_m - \Cal H^t_{n_0}| dt + \int_0^1 |\Cal H^t_{n_0}| dt
 $$
 
 $$
 \leq D(\Phi^m ,\Phi^{n_0}) + \int_0^1 |\Cal H^t_{n_0}| dt
 $$
 
 for any $n_0$. Hence if $n_0$ is sufficiently large,  $l((\Phi^n)^{-1})$  and 
 $\int_0^1 |\Cal H^t_m| dt$  are bounded. \pf

\vskip .1in

{\bf Proof of sub-lemma 2}

\vskip .1in
$$
G =\int_0^1 |\Cal H_m( \dot {\mu}_s^{n_0}) \circ \mu_s^n - \Cal H_m( \dot {\mu}_s^{n_0} )\circ \mu_s^m|ds 
$$

$$
\leq \int_0^1 |\Cal H_m( \dot {\mu}_s^{n_0}) \circ \mu_s^n - \Cal H_{m_0}( \dot {\mu}_s^{n_0} )\circ \mu_s^n|ds 
$$

$$+\int_0^1 |\Cal H_{m_0}( \dot {\mu}_s^{n_0}) \circ \mu_s^n - 
\Cal H_{m_0}( \dot {\mu}_s^{n_0} )\circ \mu_s^m|ds  
$$

$$
+\int_0^1 |\Cal H_{m_0}( \dot {\mu}_s^{n_0}) \circ \mu_s^m - 
\Cal H_m( \dot {\mu}_s^{n_0} )\circ \mu_s^m|ds 
$$

for some $m_0$.

Exactly like in the proof of sub-lemma 1
$$
G(t, n,m) \leq 2|\Cal H_m^t - \Cal H_{m_0}^t|. l(\Psi^{n_0})^{-1})) + F
$$

where

$$
F = \int_0^1 |\Cal H_{m_0}( \dot {\mu}_s^{n_0}) \circ \mu_s^n - 
\Cal H_{m_0}( \dot {\mu}_s^{n_0} )\circ \mu_s^m|ds
$$
By uniform continuity of $\Cal H_{m_0}( \dot {\mu}_s^{n_0})$, $F \to 0 $ when $n,m \to \infty$ since $\mu^n_t$ is
Cauchy.

\vskip .1in

By similar arguments as in the sub-lemma 1, $G \to 0$ and hence $L \to 0$ when $m,n \to \infty$. \pf

\vskip .1in

This concludes the proof of that $SSympeo(M,\omega)$ is a group. \pf
 
 \vskip .1in
 
The fact that it is arcwise connected in the ambiant topology of $Homeo(M)$ is obvious from the definition.

$Hameo(M,\omega)$ is a normal subgroup of $SSympeo(M,\omega)$ since it is normal in $Sympeo(M,\omega)$ [11].

Let $h,g \in SSympeo(M,\omega)$ and let  $\Phi^n ,\Psi^n$
be symplectic isotopies which form Cauchy sequences and $C^0$ converge to $h,g$. By the main lemma
 the  sequence $\Phi^n.\Psi^n.(\Phi^n)^{-1} (\Phi^n)^{-1}$ is a Cauchy sequence. It obviously converges $C^0$ to
  the commutator $hgh^{-1}g^{-1}\in SSympeo(M,\omega)$.
  
\vskip .1in

 It is a standard fact that $\Phi^n.\Psi^n.(\Phi^n)^{-1} (\Phi^n)^{-1}$ is
a hamiltonian isotopy.

Indeed let $\phi_t$ and $\psi_t$ be symplectic isotopies, and let  $\sigma_t = \phi_t \psi_t \phi_t^{-1}\psi_t^{-1}$,
 then 
 $$
 \dot {\sigma}_t = X_t + Y_t + Z_t + U_t
 $$
 
with $X_t = \dot {\phi}_t$ , $Y_t = (\phi_t)_*\dot {\psi}_t$, $Z_t = - (\phi_t \psi_t \phi_t^{-1})_*\dot {\phi}_t$,
and $U_t = - (\sigma_t)_* \dot {\psi}_t$.

By proposition 5, $i(X_t + Z_t)\omega$ and $i(Y_t + U_t)\omega$ are exact 1-forms.
Hence $\sigma_t$ is a hamiltonnian isotopy.

\vskip .1in

By proposition 1, the metric $D$ coincides with  the one for hamiltonian isotopies. Hence 
 $\Phi^n.\Psi^n.(\Phi^n)^{-1} (\Phi^n)^{-1}$ is a Cauchy sequence for $d_{ham}$. Therefore:
$[SSympeo(M,\omega), SSympeo(M,\omega)] \subset Hameo(M,\omega)]$.

This finishes the proof of the main result. \pf

\vskip .3in

{\bf Appendix}

\vskip .1in

For the convenience of the reader, we give here the proofs of propositons 3, 4, and 5.

\vskip .1in

{\bf Proof of proposition 3}

\vskip .1in

Let $\theta$ be a p-form, $X$ a vector field and $\phi$ a diffeomorphism. For any $x\in M$ and any vector
 fields $Y_1,..Y_{p-1}$, we have:
 
 $(\phi^{-1})^*[i_X\phi^*\theta](x)(Y_1, ..., Y_{p-1})  =
  (i_X\phi^*\theta)(\phi^{-1}(x))(D_x\phi^{-1}(Y_1(x),...(D_x\phi^{-1}(Y_{p-1}(x))$
  
  = $(\phi^*\theta)(\phi^{-1}(x))(X_{\phi^{-1}(x)}, D_x\phi^{-1}(Y_1(x)),...(D_x\phi^{-1}(Y_{p-1}(x))$
  
  = $\theta(\phi(\phi^{-1}(x))(D_{\phi^{-1}(x)}\phi (X_{\phi^{-1}(x)}), D_{\phi^{-1}(x)}\phi D_x\phi^{-1}(Y_1(x)),...
  D_{\phi^{-1}(x)}\phi D_x\phi^{-1}(Y_{p-1}(x)$
  
  = $\theta(x)((\phi_*X)_x, Y_1(x),.. Y_{p-1}(x))$
  
  = $(i(\phi_*X)\theta)(x)(Y_1,..,Y_{p-1})$
  
  since $D_{\phi^{-1}(x)}\phi D_x\phi^{-1} = D_x (\phi \phi^{-1}) = id$.
  
  \vskip .1in
  
  Therefore $(\phi^{-1})^*[i_X\phi^*\theta] =i(\phi_*X))\theta$
  
  \pf
  
  \vskip .1in
  
  {\bf Proof of proposition 4}
  
  \vskip .1in
  
  This is just the chain rule. See [6] page 145. \pf
  
  \vskip .1in
  
  {\bf Proof of proposition 5}
  
  \vskip .1in
   For a fixed $t$, we have
   $$
    \frac {d}{ds} \phi_s^*\theta_t = \phi_s^*( L_{\dot {\phi}_s} \theta_t,
   $$ 
   
   where $L_X$ is the Lie derivative in the direction $X$. Since $\theta$ is closed, we have:
 $$
  \frac {d}{ds} \phi_s^*\theta_t = \phi_s^*( d i_{\dot {\phi_s}}\theta_t) =
   d(\phi_s^*(\theta_t(\dot \phi_s)) = d(\theta_t(\dot \phi_s) \circ \phi_s). $$
 Hence for every $u$
 $$
 \phi_u^*\theta_t - \theta_t = \int_0^u \frac {d}{ds} \phi_s^*\theta_t  ds= d (\int_0^u (\theta_t(\dot {\phi}_s)
  \circ \phi_s) ds
 $$
 
 Now set $u = t$. \pf

{\bf References}

\vskip .2in

[1]~ A. Banyaga, {\it  On the group of symplectic homeomorphisms}

C. R. Acad. Sci. Paris Ser. 1 346(2008) 867-872

[2]~ A. Banyaga ,{\it A Hofer-like metric on the group of symplectic diffeomorphisms},

Preprint 2008.

[3]~ A. Banyaga, {\it Sur la structure du groupe des diff\'eomorphismes qui
pr\'eservent une forme symplectique},

Comment. Math. Helv. 53(1978) pp.174--227.

[4]~ A. Banyaga, {\it The structure of classical diffeomorphisms groups},

Mathematics and its applications vol 400.

Kluwer Academic Publisher's Group, Dordrecht, The Netherlands (1997).

[5]~ A. Banyaga, P. Spaeth,{\it The group of contact homeomorphisms},

Preprint, 2008.

[6]~ A. Banyaga. P. Spaeth, {\it The $C^0$ contact topology and the group of contact homeomorphisms}

Preprint, 2008

[7]~ A. Fathi, {\it Structure of the group of homeomorphisms preserving a good measure on a 
compact manifold}

Ann. Scient. Ec. Norm. Sup. 13(1980), 45-93.

[8]~ H. Hofer, E. Zehnder {\it Symplectic invariants and hamiltonian dynamics},

Birkhauser Advanced Texts, Birkhauser Verlag (1994)
 
[9]~ F. Lalonde, D. McDuff, L. Polterovich , {\it Topological rigidity of hamiltonian
 loops and quantum homology},
 
 Invent. math. 135(1999) 369-385.
 
[10]~ S. Muller, {\it The group  of hamiltonian homeomorphisms in the $L^\infty$ norm},

J. Korean Math.Soc. to appear.
 
[11]~Y-G. Oh and S. Muller {\it The group of hamiltonian homeomorphisms

and $C^0$-symplectic topology},
 
J. Symp.Geometry 5(2007) 167- 225.
 
[12]~ F. Warner {\it Foundations of differentiable manifolds and Lie groups},

Scott, Foresman and Company (1971).

 \vskip .3in
\baselineskip 12pt
\noindent Department of Mathematics \newline
\noindent The Pennsylvania State University \newline
\noindent University Park, PA 16802

\end